\documentstyle[12pt,amssymb]{article}

\begin{document}

\begin{center}
{\large {\bf Motives for elliptic modular groups}}
\end{center}

\begin{center}
By {\sc Takashi Ichikawa} 
\end{center}

\begin{center}
{\bf Abstract} 
\end{center}

In the study of the arithmetic structure of elliptic modular groups 
which are the fundamental groups of compactified modular curves, 
these truncated group algebras and their direct sums are considered 
to construct elliptic modular motives. 
Our main result is a new theory of Hecke operators on these motives 
which gives a congruence relation to the Galois action, 
and their motivic decomposition. 
Using our Hecke theory, 
we show that elliptic modular motives are the direct sums of pure motives 
over certain number fields. 
This fact implies a kind of algebraicity on iterated Shimura integrals, 
i.e., multiple $L$-values of cusp forms of weight $2$, 
and on the periods of modular Ceresa cycles. 
\vspace{2ex}

\begin{center}
{\bf 1. Introduction}
\end{center}

Under the influence of Grothendieck's consideration on motives, 
Deligne [D2] started the motivic theory of the fundamental groups 
$\pi_{1}(X)$ of algebraic varieties $X$, 
and he showed that the motives for 
$\pi_{1} \left( {\mathbb P}^{1} - \{ 0, 1, \infty \} \right)$ 
have rich structure as mixed Tate motives. 
The aim of this paper is to study {\it elliptic modular motives} 
constructed from the elliptic modular group 
$\pi_{1} \left( \overline{M}_{n}({\mathbb C}) \right)$, 
where $\overline{M}_{n}$ is the compactified modular curve 
over ${\mathbb Z} \left[ 1/n \right]$ of level $n \geq 3$. 

In order to construct motives over ${\mathbb Q}$ 
according to Deligne's theory, 
we define the Betti realization of an elliptic modular motive 
as the direct sum of the truncated group algebras of 
$\pi_{1} \left( \overline{M}_{n}({\mathbb C}) \right)$ 
whose base points run through the cusps on $\overline{M}_{n}$. 
This $l$-adic structure with action of the absolute Galois group 
${\rm Gal} \left( \overline{\mathbb Q} / {\mathbb Q} \right)$ of ${\mathbb Q}$ 
is derived from the natural ${\mathbb Z}[1/n]$-scheme structure 
on $\overline{M}_{n}$, 
and the Hodge structure is described by iterated integrals 
of $1$-forms on $\overline{M}_{n}({\mathbb C})$. 
This motive has mixed structure whose pure components are subquotients of 
$H_{1} \left( \overline{M}_{n}({\mathbb C}); {\mathbb Q} \right)^{\otimes i}$, 
and hence these components consist of tensor products of the dual motives 
attached to cusp forms of weight $2$ of level $n$. 

Our main result is to construct a new theory of Hecke operators 
which act on this elliptic modular motive and satisfy a congruence relation 
to the Galois action. 
As is well known, 
the congruence relation in the ordinary (co)homology theory is obtained by 
realizing Hecke operators as algebraic correspondences (cf. [Sm, D1]). 
However, this realization cannot work in the motivic $\pi_{1}$-case. 
Therefore, 
we define the Hecke action on the elliptic modular motive as the sum of 
operators which correspond to taking quotients of the universal elliptic curve 
by its finite subgroups. 
Then using results on algebraic fundamental groups, 
we show that this action is well-defined and 
satisfies the congruence relation. 
This Hecke theory has the following applications: 
\begin{itemize}

\item 
Our Hecke action preserves the motivic structure of 
the elliptic modular motive. 
Since iterated integrals of holomorphic $1$-forms 
on $\overline{M}_{n}({\mathbb C})$, 
known as iterated Shimura integrals of weight $2$ [M1, 2], 
are regarded as periods of this motive, 
there is Hecke action on the space of these iterated integrals. 
This gives a partial solution to a problem posed in [M2, 3.3] 
defining such action on the space of iterated Shimura integrals 
of any weight. 

\item 
The elliptic modular motive is decomposed to the direct sum of 
submotives corresponding to scalar-valued 
representations of the Hecke algebra. 
On these components, 
Hecke operators act as scalar plus nilpotent matrices, 
and hence the study of the motivic structure becomes easier. 
In particular, using the above congruence relation 
we show that this motive becomes the direct sum of pure motives 
over a number field. 

\end{itemize}
The last assertion is applied to expressing multiple $L$-values, 
represented as iterated Shimura integrals, 
of cusp forms of weight $2$ by their ordinary $L$-values, 
and to showing the algebraicity of periods integrals of divisors 
on modular curves induced from Ceresa cycles. 
\vspace{2ex}

\begin{center}
{\bf 2. Elliptic modular motives}
\end{center}

2.1. \ 
In this section, 
we construct motives for elliptic modular groups 
based on Deligne's theory [D1]. 
Let $\pi_{1}$ be the topological fundamental group 
$\pi_{1} \left( R; x \right)$ of a Riemann surface $R$ with base point $x$. 
Then $\pi_{1}$ is finitely generated. 
Denote by ${\mathbb Q} \left[ \pi_{1} \right]$ 
the group ${\mathbb Q}$-algebra of $\pi_{1}$ with the augmentation ideal $I$, 
and by 
$$
{\mathbb Q} \left[ \pi_{1} \right]^{\wedge} = 
\lim_{\longleftarrow} {\mathbb Q} \left[ \pi_{1} \right] / I^{m}, 
$$
the completion of ${\mathbb Q} \left[ \pi_{1} \right]$ 
with the augmentation ideal $\widehat{I}$ given as the completion of $I$. 
Let 
$$
\Delta : {\mathbb Q} \left[ \pi_{1} \right]^{\wedge} \rightarrow 
{\mathbb Q} \left[ \pi_{1} \right]^{\wedge} \widehat{\otimes} \ 
{\mathbb Q} \left[ \pi_{1} \right]^{\wedge}  
$$
be the diagonal homomorphism which is a continuous algebra homomorphism 
induced from $g \mapsto g \otimes g$ $(g \in \pi_{1})$. 
Then the Malcev Lie algebra ${\rm Lie}(\pi_{1})$ of $\pi_{1}$ 
is defined as the set 
$$
\left\{ \left. a \in {\mathbb Q} \left[ \pi_{1} \right]^{\wedge} \ \right| \ 
\Delta(a) = 1 \otimes a + a \otimes 1 \right\} 
$$
of primitive elements in which the bracket product $[a, b]$ 
is given by $ab - ba$. 
Then the exponential map 
$$
\exp(a) = \sum_{i=0}^{\infty} \frac{a^{i}}{i!} 
$$
gives a bijection from ${\rm Lie}(\pi_{1})$ onto the set 
$$
\left\{ g \in {\mathbb Q} \left[ \pi_{1} \right]^{\wedge} \ \left| \ 
g - 1 \in \widehat{I}, \ \Delta(g) = g \otimes g \right. \right\} 
$$
of grouplike elements, 
and its inverse map is the logarithmic map 
$$
\log(g) = - \sum_{i=1}^{\infty} \frac{(1-g)^{i}}{i}. 
$$ 

For a positive integer $N$, 
let 
$$
A_{N}(R; x) = {\mathbb Q} \left[ \pi_{1}(R; x) \right] \left/ I^{N+1} \right. 
$$
be the $N$th truncated group algebra of $\pi_{1}(R; x)$ over ${\mathbb Q}$. 
Then the $N$th truncated Malcev Lie algebra 
${\rm Lie}_{N} (\pi_{1}) = 
{\rm Lie}_{N} \left( \pi_{1} \left( R; x \right) \right)$ 
is defined as the image of ${\rm Lie} \left( \pi_{1}(R; x) \right)$ 
into $A_{N}(R; x)$ by the natural projection. 
The corresponding unipotent algebraic group 
$G_{N} = G_{N} \left( \pi_{1} \left( R; x \right) \right)$ 
over ${\mathbb Q}$ is characterized by that $G_{N} \left( {\mathbb Q} \right)$ 
consists of the images of grouplike elements into $A_{N}(R; x)$ 
by the natural projection. 
Further, for any field $K$ of characteristic $0$, 
the exponential and logarithmic maps give bijections 
${\rm Lie}_{N} (\pi_{1}) \otimes K \stackrel{\sim}{\leftrightarrow} 
G_{N} \left( K \right)$ 
which are inverse to each other. 

For a prime $l$, 
denote by $\pi_{1}^{l}$ the $l$-adic (algebraic) fundamental group. 
For results on this subject, see [G]. 
Then for $x \in R$, 
the $l$-adic completion of $\pi_{1}(R; x)$ is canonically isomorphic to 
$\pi_{1}^{l}(R; x)$, 
and hence the natural group homomorphism 
$$
\pi_{R} : \pi_{1}(R; x) \rightarrow G_{N} \left( {\mathbb Q} \right) 
\subset A_{N}(R; x) 
$$
gives rise to a group homomorphism 
$$
\pi_{1}^{l}(R; x) \rightarrow G_{N} \left( {\mathbb Q}_{l} \right) 
\subset A_{N}(R; x) \otimes {\mathbb Q}_{l} 
$$
which we denote by the same symbol. 

Let $\phi : R' \rightarrow R$ be a morphism between Riemann surfaces 
sending $x' \in R'$ to $x \in R$. 
Then we have the associated group homomorphism 
$$
\pi_{1} (\phi) : \pi_{1} \left( R'; x' \right) \rightarrow \pi_{1} (R; x) 
$$
which gives rise to a ${\mathbb Q}$-algebra homomorphism 
$$
A_{N} (\phi) : A_{N} \left( R'; x' \right) \rightarrow A_{N} (R; x)  
$$
and a group homomorphism 
$$
\pi_{1}^{l} (\phi) : \pi_{1}^{l} \left( R'; x' \right) \rightarrow 
\pi_{1}^{l} (R; x). 
$$

2.2. \ 
Let $n$ be an integer $\geq 3$, 
denote by $P_{n}$ the set of primitive $n$th roots of $1$, 
and by $R_{n}$ the subring ${\mathbb Z}[ 1/n, \zeta ]$ of ${\mathbb C}$, 
where $\zeta$ is an element of $P_{n}$. 
Let $M_{n}$ be the affine modular curve over ${\mathbb Z}[1/n]$ 
which classifies elliptic curves $E$ with level $n$ structure, 
i.e., isomorphism 
$$
\lambda : \left( {\mathbb Z}/n {\mathbb Z} \right)^{\oplus 2} 
\stackrel{\sim}{\rightarrow} 
E[n] \stackrel{\rm def}{=} {\rm Ker}(n : E \rightarrow E). 
$$
For each $\zeta \in P_{n}$, 
let $M_{\zeta}$ be the modular curve over $R_{n}$ 
which classifies $(E, \lambda)$ such that 
$e(a, b) = \zeta^{\psi(\lambda^{-1}(a), \lambda^{-1}(b))}$, 
where $e$ denotes the Weil pairing, 
and $\psi$ denotes the standard symplectic form. 
Then each $M_{\zeta}$ is geometrically irreducible, 
and 
$$
M_{n} \otimes R_{n} = \bigsqcup_{\zeta \in P_{n}} M_{\zeta}. 
$$
Further, 
$M_{\zeta}({\mathbb C}) \cong H / \Gamma(n)$, 
where $H = \{ \tau \in {\mathbb C} \ | \ {\rm Im}(\tau) > 0 \}$ 
denotes the Poincar\'{e} upper half plane, 
and 
$\Gamma(n) = {\rm Ker} \left( SL_{2}({\mathbb Z}) \rightarrow 
SL_{2} \left( {\mathbb Z} / n {\mathbb Z} \right) \right)$ 
denotes the principal congruence subgroup of level $n$. 
Let $\overline{M}_{\zeta} = M_{\zeta} \cup C_{\zeta}$ 
be the compactified modular curve over $R_{n}$, 
where $C_{\zeta}$ denotes the scheme of cusps which is finite and \'{e}tale 
over $R_{n}$ and consists of $R_{n}$-rational points (see [I]). 
Then $\overline{M}_{\zeta}({\mathbb C}) \cong \overline{H} / \Gamma(n)$, 
where $\overline{H} = H \cup {\mathbb P}^{1}({\mathbb Q})$, 
and there is a unique proper smooth curve $\overline{M}_{n}$ 
over ${\mathbb Z}[ 1/n ]$ such that 
$\overline{M}_{n} \otimes R_{n} = 
\bigsqcup_{\zeta \in P_{n}} \overline{M}_{\zeta}$. 

In what follows, fix a positive integer $N$. 
Then the Betti realization $V = V_{N}$ 
of the {\it $N$th truncated elliptic modular motive of level $n$} 
is defined as the direct sum of the $N$th truncated group algebra of 
$\pi_{1} \left( \overline{M}_{\zeta}({\mathbb C}); c \right)$, 
where $\zeta$ and $c$ run through $P_{n}$ and $C_{\zeta}$ respectively: 
$$
V = V_{N} = \bigoplus_{\zeta, c} 
A_{N} \left( \overline{M}_{\zeta}({\mathbb C}); c \right). 
$$ 
We regard $V$ together with the following structure as a motive. 

First, we consider the $l$-adic structure on $V$, i.e., 
the ${\rm Gal} \left( \overline{\mathbb Q}/{\mathbb Q} \right)$-module 
structure on $V \otimes {\mathbb Q}_{l}$. 
For $x \in \overline{M}_{\zeta} \left( \overline{\mathbb Q} \right)$, 
$\pi_{1}^{l} \left( \overline{M}_{\zeta}({\mathbb C}); x \right) \cong 
\pi_{1}^{l} \left( \overline{M}_{\zeta} \otimes \overline{\mathbb Q}; 
x \right)$ 
is the automorphism group of the fiber systems 
$\left\{ f_{Y}^{-1}(x) \right\}$ for Galois coverings 
$f_{Y} : Y \rightarrow \overline{M}_{\zeta} \otimes \overline{\mathbb Q}$ 
of $l$-power degree. 
Therefore, 
for $\sigma \in {\rm Gal} \left( \overline{\mathbb Q}/{\mathbb Q} \right)$, 
the correspondence 
$$
\gamma \mapsto \sigma \circ \gamma \circ \sigma^{-1} \ 
\left( \gamma \in \pi_{1}^{l} \left( \overline{M}_{\zeta}({\mathbb C}); 
x \right) \right)
$$
together with $\pi_{\overline{M}_{\zeta}({\mathbb C})}$ 
give a ${\mathbb Q}_{l}$-algebra homomorphism 
$$
\eta(\sigma) : 
A_{N} \left( \overline{M}_{\zeta}({\mathbb C}); x \right) \otimes 
{\mathbb Q}_{l} 
\rightarrow 
A_{N} \left( \overline{M}_{\sigma(\zeta)}({\mathbb C}); \sigma(x) \right) 
\otimes {\mathbb Q}_{l}. 
$$
Since $\bigsqcup_{\zeta} C_{\zeta}$ is stable under the Galois action 
over ${\mathbb Q}$, 
${\rm Gal} \left( \overline{\mathbb Q}/{\mathbb Q} \right)$ acts naturally on 
$V \otimes {\mathbb Q}_{l}$. 

Second, following Hain [H1, 2] (see also [Mo]), 
we consider weight and Hodge filtrations of $V$. 
By results of Chen [Ch1, 2], for each $c \in C_{\zeta}$, 
$A_{N} \left( \overline{M}_{\zeta}({\mathbb C}); c \right) \otimes 
{\mathbb C}$ 
is dual to the ${\mathbb C}$-vector space 
$B_{N} \left( \overline{M}_{\zeta}({\mathbb C}); c \right)$ 
consisting of finite sums of iterated integrals 
$$
\int w_{1} \cdots w_{r} \ 
\left( \mbox{$w_{i} :$ smooth $1$-forms on 
$\overline{M}_{\zeta}({\mathbb C})$, $r \leq N$} \right) 
$$
which are homotopy functional on 
$\left\{ \mbox{loops based at $c$} \right\}$. 
Hence $V \otimes {\mathbb C}$ is dual to 
$$
B_{N} = \bigoplus_{\zeta, c} 
B_{N} \left( \overline{M}_{\zeta}({\mathbb C}); c \right), 
$$
where $\zeta, c$ run through $P_{n}, C_{\zeta}$ respectively. 
Hain [H1, 2] defined the weight and Hodge filtrations of $B_{N}$ as 
\begin{eqnarray*}
W^{l} \left( B_{N} \right) 
& = & 
\left\{ \begin{array}{ll} 
B_{l} & (l \leq N), 
\\
B_{N} & (l \geq N), 
\end{array} \right. 
\\
F^{k} \left( B_{N} \right) 
& = & 
\bigoplus_{\zeta, c} 
\left\{ \begin{array}{l}
\mbox{sums of iterated integrals of $1$-forms with $\geq k$ $dz$'s,} 
\\
\mbox{where $z$ is a holomorphic coordinate on $\overline{M}_{\zeta}$} 
\end{array} \right\}, 
\end{eqnarray*}
and hence by the duality, 
$V$ has weight filtration and $V \otimes {\mathbb C}$ has Hodge filtration. 
Furthermore, 
Hain defined de Rham structure on $V$ which is induced from 
the de Rham cohomology group 
$H_{\rm DR}^{1} \left( \overline{M}_{\zeta}; {\mathbb Q}(\zeta) \right)$ 
and is compatible with the Hodge filtration on $V \otimes {\mathbb C}$. 

Finally, we only notice that the crystalline structure on $V$ 
is induced from that on the first cohomology groups of 
$\overline{M}_{\zeta}$ $(\zeta \in P_{n})$ 
and is constructed by Shiho [S]. 
\vspace{2ex}

\begin{center}
{\bf 3. Hecke operators and congruence relation}
\end{center}

3.1. \ 
Denote by $F_{\zeta}$ the function field of $M_{\zeta}$, 
and by $\left( {\cal E}_{\zeta}, \lambda_{\zeta} \right)$ the universal 
elliptic curve with level $n$ structure over $M_{\zeta}$. 
Let $G$ be a subgroup of ${\cal E}_{\zeta} \otimes \overline{F}_{\zeta}$ 
such that its order $|G|$ is prime to $n$. 
Take a finite extension $K$ of $F_{\zeta} \cdot \overline{\mathbb Q}$ 
such that $G$ is defined over $K$, 
and let $\overline{X}_{\zeta, G}$ be a proper smooth curve over 
$\overline{\mathbb Q}$ whose function field is $K$. 
Denote by $\phi_{\zeta, G} : \overline{X}_{\zeta, G} \rightarrow 
\overline{M}_{\zeta} \otimes \overline{\mathbb Q}$ the natural morphism, 
and put 
$X_{\zeta, G} = 
\phi_{\zeta, G}^{-1} \left( M_{\zeta} \otimes \overline{\mathbb Q} \right)$. 
\vspace{2ex}

{\sc Proposition 3.1.} 
\begin{it}
For any $z \in \overline{X}_{\zeta, G}({\mathbb C})$, 
the ${\mathbb Q}$-algebra homomorphism 
$$
A_{N} \left( \phi_{\zeta, G} \right) : 
A_{N} \left( \overline{X}_{\zeta, G}({\mathbb C}); z \right) \rightarrow 
A_{N} \left( \overline{M}_{\zeta}({\mathbb C}); \phi_{\zeta, G}(z) \right)
$$
is surjective. 
\end{it} 
\vspace{2ex}

{\it Proof.} 
Take $\gamma \in \pi_{1} \left( \overline{M}_{\zeta}({\mathbb C}); 
\phi_{\zeta, G}(z) \right)$. 
Since 
$\phi_{\zeta, G} : \overline{X}_{\zeta, G} \rightarrow 
\overline{M}_{\zeta} \otimes \overline{\mathbb Q}$ 
is a finite morphism, 
${\rm Im} \left( \pi_{1} \left( \phi_{\zeta, G} \right) \right)$ 
is finite index in 
$\pi_{1} \left( \overline{M}_{\zeta}({\mathbb C}); 
\phi_{\zeta, G}(z) \right)$. 
Therefore, 
there are a positive integer $t$ and 
$\widetilde{\gamma} \in \pi_{1} \left( \overline{X}_{\zeta, G}({\mathbb C}); 
z \right)$ 
such that 
$\gamma^{t} = \pi_{1} \left( \phi_{\zeta, G} \right) 
\left( \widetilde{\gamma} \right)$. 
Since $A_{N} \left( \phi_{\zeta, G} \right)$ is 
a ${\mathbb Q}$-algebra homomorphism, 
\begin{eqnarray*}
\gamma 
& = &  
\exp \left( \frac{1}{t} \log \left( \pi_{\overline{M}_{\zeta}({\mathbb C})} 
\left( \pi_{1} \left( \phi_{\zeta, G} \right) 
\left( \widetilde{\gamma} \right) \right) \right) \right) 
\\ 
& = & 
A_{N} \left( \phi_{\zeta, G} \right) \left( \exp 
\left( \frac{1}{t} \log \left( \pi_{\overline{X}_{\zeta, G}({\mathbb C})} 
\left( \widetilde{\gamma} \right) \right) \right) \right) 
\end{eqnarray*}
belongs to the image of $A_{N} \left( \phi_{\zeta, G} \right)$, 
and hence this map is surjective. \ $\square$ 
\vspace{2ex}

The quotient 
$\left( {\cal E}_{\zeta}, \lambda_{\zeta} \right) \otimes K / G$ 
is an elliptic curve over $K$ with level $n$ structure which gives rise to 
a $K$-rational point on $M_{\zeta^{|G|}}$, 
and hence we have the associated morphism 
$\psi_{\zeta, G} : \overline{X}_{\zeta, G} \rightarrow 
\overline{M}_{\zeta^{|G|}} \otimes \overline{\mathbb Q}$. 
\vspace{2ex}

{\sc Theorem 3.2.} 
\begin{it}
Assume that $|G|$ is a prime $p$ not dividing $nl$, 
and let $z$ be a point on 
$\overline{X}_{\zeta, G} \left( \overline{\mathbb Q} \right)$. 
\begin{itemize}

\item[{\rm (1)}] 
There exists a group homomorphism 
$$
\tau_{\zeta, G, z}^{l} : 
\pi_{1}^{l} \left( \overline{M}_{\zeta} \otimes \overline{\mathbb Q}; 
\phi_{\zeta, G}(z) \right) 
\rightarrow 
\pi_{1}^{l} \left( \overline{M}_{\zeta^{p}} \otimes \overline{\mathbb Q}; 
\psi_{\zeta, G}(z) \right) 
$$
satisfying that 
$\tau_{\zeta, G, z}^{l} \circ \pi_{1}^{l} \left( \phi_{\zeta, G} \right) 
= \pi_{1}^{l} \left( \psi_{\zeta, G} \right)$. 
Actually, $\tau_{\zeta, G, z}^{l}$ is induced from either 
the Frobenius or Vershiebung at an extension of $p$. 

\item[{\rm (2)}] 
If $\phi_{\zeta, G}(z) \in C_{\zeta}$, 
then for any 
$z' \in \overline{X}_{\zeta, G} \left( \overline{\mathbb Q} \right)$ such that 
$\phi_{\zeta, G}(z') = \phi_{\zeta, G}(z)$, 
$\psi_{\zeta, G}(z') = \psi_{\zeta, G}(z)$. 
\end{itemize}
\end{it} 

{\it Proof.} 
Let $\overline{p}$ be an extension of $p$ as a place of 
$\overline{F}_{\zeta}$,  
and $\overline{k}$ be the residue field of the valuation ring $\overline{R}$ 
of $\overline{F}_{\zeta}$ at $\overline{p}$ such that 
${\cal E}_{\zeta} \otimes \overline{k}$ corresponds to 
the generic point on $\overline{M}_{\zeta} \otimes k$, 
where $R$ and $k$ denote the valuation ring of $\overline{\mathbb Q}$ 
at $\overline{p}$ and its residue field respectively. 
Let $\overline{Y}_{\zeta, G}$ be the normalization of 
$\overline{M}_{\zeta} \otimes R$ in $K$, 
and 
$\phi_{\zeta, G} : \overline{Y}_{\zeta, G} \rightarrow 
\overline{M}_{\zeta} \otimes R$ 
be the natural morphism. 
Since $\overline{Y}_{\zeta, G}$ is proper over $R$, 
and $\overline{M}_{\zeta} \otimes R$ is proper smooth over $R$, 
we have the following commutative diagram: 
$$
\begin{array}{ccccc} 
\pi_{1}^{l} \left( \overline{Y}_{\zeta, G} \otimes \overline{\mathbb Q}; 
z \right) 
& \rightarrow & 
\pi_{1}^{l} \left( \overline{Y}_{\zeta, G}; z \right) 
& \stackrel{\sim}{\leftarrow} & 
\pi_{1}^{l} \left( \overline{Y}_{\zeta, G} \otimes k; z \right) 
\\
\downarrow & & \downarrow & & \downarrow 
\\
\pi_{1}^{l} \left( \overline{M}_{\zeta} \otimes \overline{\mathbb Q}; 
\phi_{\zeta, G}(z) \right) 
& \stackrel{\sim}{\rightarrow} & 
\pi_{1}^{l} \left( \overline{M}_{\zeta} \otimes R; \phi_{\zeta, G}(z) \right) 
& \stackrel{\sim}{\leftarrow} & 
\pi_{1}^{l} \left( \overline{M}_{\zeta} \otimes k; \phi_{\zeta, G}(z) \right). 
\end{array}
$$ 

Let $Y_{\zeta, G}$ be the inverse image of $M_{\zeta} \otimes R$ 
by $\phi_{\zeta, G}$ which is an affine normal scheme over $R$ 
such that $Y_{\zeta, G} \otimes \overline{\mathbb Q} \cong X_{\zeta, G}$. 
Since ${\cal E}_{\zeta} \otimes K / G$ has good reduction 
over any discrete valuation ring of $K$ containing $A$, 
where ${\rm Spec}(A)$ is an open subset of $Y_{\zeta, G}$, 
$\psi_{\zeta, G}|_{X_{\zeta, G}}$ is uniquely extended to a morphism 
$Y_{\zeta, G} \rightarrow M_{\zeta^{|G|}} \otimes R$ 
which we denote by $\psi_{\zeta, G}$. 
Further, by the theory of Tate curves, 
the Neron model of 
$\psi_{\zeta, G}^{*} \left( {\cal E}_{\zeta^{|G|}} \right)$ over 
$\overline{Y}_{\zeta, G}$ gives rise to a morphism 
$\overline{Y}_{\zeta, G} \rightarrow \overline{M}_{\zeta^{|G|}} \otimes R$ 
as a unique extension of the above $\psi_{\zeta, G}$ 
which we denote by $\psi_{\zeta, G}$ also. 

Let $G_{\overline{p}}$ be the reduction modulo $\overline{p}$ of the model 
of $G$ as a subgroup scheme of 
${\cal E}_{\zeta} \otimes \overline{R}$ flat over $\overline{R}$. 
Then $G_{\overline{p}}$ is a subgroup scheme of the ordinary elliptic curve 
${\cal E}_{\zeta} \otimes \overline{k}$, 
and hence one of the following two cases necessarily happens: 
\vspace{2ex}

{\it Case} 1. $G_{\overline{p}}$ is connected. 
In this case, 
$G_{\overline{p}}$ is the kernel of the Frobenius homomorphism 
${\cal E}_{\zeta} \otimes \overline{k} \rightarrow 
\left( {\cal E}_{\zeta} \otimes \overline{k} \right)^{(p)}$, 
and hence 
$$
\psi_{\zeta, G}^{*} \left( \left( {\cal E}_{\zeta^{p}}, 
\lambda_{\zeta^{p}} \right) \otimes \overline{k} \right) 
\cong 
\phi_{\zeta, G}^{*} \left( \left( {\cal E}_{\zeta}, \lambda_{\zeta} \right) 
\otimes \overline{k} \right)^{(p)}. 
$$
Let $\left( {\cal G}_{\zeta}, \lambda_{\zeta} \right)$ 
be the generalized elliptic curve with level $n$ structure 
over $\overline{M}_{\zeta}$ (cf. [DR, KM]). 
Then we have 
$$
\psi_{\zeta, G}^{*} \left( \left( {\cal G}_{\zeta^{p}}, 
\lambda_{\zeta^{p}} \right) \otimes \overline{k} \right) 
\cong 
\phi_{\zeta, G}^{*} \left( \left( {\cal G}_{\zeta}, \lambda_{\zeta} \right) 
\otimes \overline{k} \right)^{(p)}, 
$$
and hence the Frobenius morphism 
$\varphi_{\overline{p}} : \overline{M}_{\zeta} \otimes k \rightarrow 
\overline{M}_{\zeta^{p}} \otimes k$ 
at $\overline{p}$ sending 
$\left( {\cal G}_{\zeta}, \lambda_{\zeta} \right) \otimes \overline{k}$ 
to 
$\left( \left( {\cal G}_{\zeta}, \lambda_{\zeta} \right) \otimes \overline{k} 
\right)^{(p)}$ 
satisfies that 
$\varphi_{\overline{p}} \circ \phi_{\zeta, G} = \psi_{\zeta, G}$ 
on $\overline{Y}_{\zeta, G} \otimes k$. 
Therefore, by the above diagram, 
$\varphi_{\overline{p}}$ gives rise to $\tau_{\zeta, G, z}^{l}$ 
satisfying the condition in (1). 
Since 
$$
\varphi_{\overline{p}} \left( \phi_{\zeta, G}(z) \ 
{\rm mod} \left( \overline{p} \right) \right) = 
\psi_{\zeta, G}(z) \ {\rm mod} \left( \overline{p} \right) 
$$ 
and the reduction map modulo $\overline{p}$ on $C_{\zeta}$ is bijective, 
the assertion (2) holds. 
\vspace{2ex}

{\it Case} 2. $G_{\overline{p}}$ is \'{e}tale over $\overline{k}$. 
In this case, 
$G_{\overline{p}}$ is the kernel of the Vershiebung map 
${\cal E}_{\zeta} \otimes \overline{k} \rightarrow 
\left( {\cal E}_{\zeta} \otimes \overline{k} \right)^{\left( p^{-1} \right)}$, 
and hence 
$$
\phi_{\zeta, G}^{*} \left( \left( {\cal G}_{\zeta}, \lambda_{\zeta} \right) 
\otimes \overline{k} \right) 
\cong 
\psi_{\zeta, G}^{*} \left( \left( {\cal G}_{\zeta^{p}}, 
\lambda_{\zeta^{p}} \right) \otimes \overline{k} \right)^{(p)}. 
$$
Then as in the Case 1, 
the Vershiebung at $\overline{p}$ gives $\tau_{\zeta, G, z}^{l}$ 
satisfying the condition in (1) and sends $\phi_{\zeta, G}(z) \in C_{\zeta}$ 
to $\psi_{\zeta, G}(z)$. 
Therefore, (2) holds also. 
\ $\square$ 
\vspace{2ex} 

{\sc Theorem 3.3.} 
\begin{it} 
Assume that $|G|$ is prime to $n$, 
and let $z$ be a point on 
$\overline{X}_{\zeta, G} \left( \overline{\mathbb Q} \right)$. 
\begin{itemize}

\item[{\rm (1)}] 
If $|G|$ is prime to $l$, 
then there exists a group homomorphism 
$$
\tau_{\zeta, G, z}^{l} : 
\pi_{1}^{l} \left( \overline{M}_{\zeta} \otimes \overline{\mathbb Q}; 
\phi_{\zeta, G}(z) \right) 
\rightarrow 
\pi_{1}^{l} \left( \overline{M}_{\zeta^{|G|}} \otimes \overline{\mathbb Q}; 
\psi_{\zeta, G}(z) \right) 
$$
satisfying that 
$\tau_{\zeta, G, z}^{l} \circ \pi_{1}^{l} \left( \phi_{\zeta, G} \right) 
= \pi_{1}^{l} \left( \psi_{\zeta, G} \right)$. 
Furthermore, 
if $\phi_{\zeta, G}(z) \in C_{\zeta}$, 
then for any 
$z' \in \overline{X}_{\zeta, G} \left( \overline{\mathbb Q} \right)$ such that 
$\phi_{\zeta, G}(z') = \phi_{\zeta, G}(z)$, 
$\psi_{\zeta, G}(z') = \psi_{\zeta, G}(z)$. 

\item[{\rm (2)}] 
There exists a unique ${\mathbb Q}$-linear map 
$$
\tau_{\zeta, G, z} : 
A_{N} \left( \overline{M}_{\zeta}({\mathbb C}); \phi_{\zeta, G}(z) \right) 
\rightarrow 
A_{N} \left( \overline{M}_{\zeta^{|G|}}({\mathbb C}); 
\psi_{\zeta, G}(z) \right) 
$$
satisfying that 
$\tau_{\zeta, G, z} \circ A_{N} \left( \phi_{\zeta, G} \right) 
= A_{N} \left( \psi_{\zeta, G} \right)$. 
Furthermore, $\tau_{\zeta, G, z}$ is independent of the choice of 
$\overline{X}_{\zeta, G}$, 
and $\tau_{\zeta, G, z} \otimes {\mathbb Q}_{l}$ is given by 
$\tau_{\zeta, G, z}^{l}$. 

\item[{\rm (3)}] 
For $\sigma \in {\rm Gal} \left( \overline{\mathbb Q} / {\mathbb Q} \right)$, 
let 
$\widetilde{\sigma} \in 
{\rm Aut} \left( \overline{F}_{\zeta} / {\mathbb Q} \right)$ 
be an extension of $\sigma$. 
Then 
$$
\tau_{\sigma(\zeta), \widetilde{\sigma}(G), \sigma(z)} \circ \eta(\sigma) 
= \eta(\sigma) \circ \tau_{\zeta, G, z} 
$$
holds on 
$A_{N} \left( \overline{M}_{\zeta}({\mathbb C}); \phi_{\zeta, G}(z) \right) 
\otimes {\mathbb Q}_{l}$. 

\end{itemize}
\end{it} 

{\it Proof.} 
The assertion (1) follows from Theorem 3.2. 
The uniqueness of $\tau_{\zeta, G, z}$ in (2) follows from Proposition 3.1, 
and we show the existence. 
The group homomorphism $\tau_{\zeta, G, z}^{l}$ gives 
a ${\mathbb Q}_{l}$-linear map 
$$
A_{N} \left( \overline{M}_{\zeta}({\mathbb C}); \phi_{\zeta, G}(z) \right) 
\otimes {\mathbb Q}_{l} 
\rightarrow 
A_{N} \left( \overline{M}_{\zeta^{|G|}}({\mathbb C}); 
\psi_{\zeta, G}(z) \right) \otimes {\mathbb Q}_{l} 
$$    
which we denote by the same symbol. 
Since 
$\tau_{\zeta, G, z}^{l} \circ A_{N} \left( \phi_{\zeta, G} \right) 
= A_{N} \left( \psi_{\zeta, G} \right)$ 
on 
$A_{N} \left( \overline{X}_{\zeta, G}({\mathbb C}); z \right) 
\otimes {\mathbb Q}_{l}$, 
and $A_{N} \left( \phi_{\zeta, G} \right)$ is 
a surjective ${\mathbb Q}$-linear map, 
the restriction $\tau_{\zeta, G, z}$ of $\tau_{\zeta, G, z}^{l}$ to 
$A_{N} \left( \overline{M}_{\zeta}({\mathbb C}); \phi_{\zeta, G}(z) \right)$ 
satisfies the condition in (2). 
The remaining assertions in (2) are clear. 
Finally, we prove (3). 
Since 
$\sigma \left( \phi_{\zeta, G} \right) = 
\phi_{\sigma(\zeta), \widetilde{\sigma}(G)}$ and 
$\sigma \left( \psi_{\zeta, G} \right) = 
\psi_{\sigma(\zeta), \widetilde{\sigma}(G)}$, 
the homomorphism 
$\eta(\sigma) \circ \tau_{\zeta, G, z}^{l} \circ \eta(\sigma)^{-1}$ 
defined on 
$\pi_{1}^{l} \left( \overline{M}_{\sigma(\zeta)} \otimes \overline{\mathbb Q}; 
\sigma \left( \phi_{\zeta, G}(z) \right) \right)$ 
satisfies that 
$$
\left( \eta(\sigma) \circ \tau_{\zeta, G, z}^{l} \circ \eta(\sigma)^{-1} 
\right) \circ \pi_{1}^{l} \left( \sigma \left( \phi_{\zeta, G} \right) \right) 
= \pi_{1}^{l} \left( \sigma \left( \psi_{\zeta, G} \right) \right). 
$$
Therefore, by (2), we have 
$\tau_{\sigma(\zeta), \widetilde{\sigma}(G), \sigma(z)} \circ \eta(\sigma) 
= \eta(\sigma) \circ \tau_{\zeta, G, z}$. 
\ $\square$ 
\vspace{2ex} 

{\sc Proposition 3.4.} 
\begin{it}
Let $G$ and $H$ be subgroups of 
${\cal E}_{\zeta} \otimes \overline{F}_{\zeta}$ and 
${\cal E}_{\zeta} \otimes \overline{F}_{\zeta} / G$ respectively 
such that $|G| \cdot |H|$ is prime to $n$. 
Denote by $\overline{H}$ the inverse image of $H$ by the projection 
${\cal E}_{\zeta} \otimes \overline{F}_{\zeta} \rightarrow 
{\cal E}_{\zeta} \otimes \overline{F}_{\zeta} / G$, 
and let $\overline{X}_{\zeta, \overline{H}}$ be a proper smooth curve 
over which $G$ and $\overline{H}$ are defined. 
Then there exist morphisms 
$\pi_{1} : \overline{X}_{\zeta, \overline{H}} \rightarrow 
\overline{X}_{\zeta, G}$ 
and 
$\pi_{2} : \overline{X}_{\zeta, \overline{H}} \rightarrow 
\overline{X}_{\zeta^{|G|}, H}$ 
such that 
$\phi_{\zeta, \overline{H}} = \phi_{\zeta, G} \circ \pi_{1}$, 
$\psi_{\zeta, \overline{H}} = \psi_{\zeta^{|G|}, H} \circ \pi_{2}$ 
and that 
$$
\tau_{\zeta^{|G|}, H, \left( \psi_{\zeta, G} \circ \pi_{1} \right)(z)} 
\circ \tau_{\zeta, G, z} 
= \tau_{\zeta, \overline{H}, z}. 
$$
\end{it}
\vspace{-2ex}

{\it Proof.}
This follows from the construction of $\phi_{*, *}$, $\psi_{*, *}$ 
and Theorem 3.3. 
\ $\square$ 
\vspace{2ex}

3.2. \ 
By Theorem 3.3, 
there exists a map $\iota_{\zeta, G} : C_{\zeta} \rightarrow C_{\zeta^{|G|}}$ 
for each $\zeta, G$ as above such that if $c \in C_{\zeta}$, 
then $\tau_{\zeta, G, *}$ gives a ${\mathbb Q}$-linear map 
$$
\tau_{\zeta, G} : A_{N} \left( \overline{M}_{\zeta}({\mathbb C}); c \right) 
\rightarrow 
A_{N} \left( \overline{M}_{\zeta^{|G|}}({\mathbb C}); 
\iota_{\zeta, G}(c) \right). 
$$
Therefore, 
for each pair $(m_{1}, m_{2})$ of positive integers such that 
$m_{1} | m_{2}$ and that $m_{1} m_{2}$ is prime to $n$, 
one can define a ${\mathbb Q}$-linear endomorphism $T(m_{1}, m_{2})$ on 
$$
V = \bigoplus_{\zeta, c} 
A_{N} \left( \overline{M}_{\zeta}({\mathbb C}); c \right)
$$ 
as the sum of $\tau_{\zeta, G}$, 
where $G$ runs through all the subgroups of 
${\cal E}_{\zeta} \otimes \overline{F}_{\zeta}$ which are of type 
$(m_{1}, m_{2})$, i.e., isomorphic to 
${\mathbb Z}/m_{1} \oplus {\mathbb Z}/m_{2}$.  
In particular, put $T(p) = T(1, p)$ for each prime $p$ not dividing $n$. 
\vspace{2ex}

{\sc Theorem 3.5.} 
\begin{it}
\begin{itemize}

\item[{\rm (1)}] 
For each prime $p$ not dividing $n$, 
$$
T(p) = F_{p} + p \left( F_{p}^{-1} \circ T(p, p) \right) 
$$ 
holds on $V \otimes {\mathbb Q}_{l}$, 
where $F_{p} \in {\rm Gal} \left( \overline{\mathbb Q}/{\mathbb Q} \right)$ 
denotes an isomorphism given by the Frobenius morphism. 

\item[{\rm (2)}] 
For any $\sigma \in 
{\rm Gal} \left( \overline{\mathbb Q}/{\mathbb Q} \right)$, 
$$
\eta(\sigma) \circ T(m_{1}, m_{2}) = T(m_{1}, m_{2}) \circ \eta(\sigma)
$$ 
holds on $V \otimes {\mathbb Q}_{l}$. 

\end{itemize}
\end{it} 

{\it Proof.} 
The assertion (1) follows from Theorem 3.2 since there are $p+1$ subgroups of 
${\cal E}_{\zeta} \otimes \overline{F}_{\zeta}$ and their reductions modulo 
$\overline{p}$ consist of the connected subgroup scheme and 
$p$ copies of the \'{e}tale subgroup scheme of 
${\cal E}_{\zeta} \otimes \overline{k}$. 
Further, (2) follows from Theorem 3.3 (3) 
since if $G$ runs through all the subgroups of 
${\cal E}_{\zeta} \otimes \overline{F}_{\zeta}$ of type $(m_{1}, m_{2})$, 
then $\widetilde{\sigma}(G)$ does so. 
\ $\square$. 
\vspace{2ex}

\begin{center}
{\bf 4. Decomposition by Hecke action}
\end{center}

4.1. \ 
In this section, 
we introduce a theory of Hecke operators which gives a motivic decomposition 
of the elliptic modular motive with congruence relation. 
Let ${\cal H}$ be the Hecke algebra which is defined as 
a ${\mathbb Q}$-algebra generated by double coset classes 
$$
\widetilde{T}(m_{1}, m_{2}) = \Gamma(n) 
\left( \begin{array}{cc} m_{1} & 0 \\ 0 & m_{2} \end{array} \right) \Gamma(n) 
$$
for positive integers $m_{1}, m_{2}$ prime to $n$ such that $m_{1} | m_{2}$. 
Then it is shown in [Sm] that ${\cal H}$ becomes a commutative ring. 
\vspace{2ex}

{\sc Proposition 4.1.} 
\begin{it} 
There exists a unique ${\mathbb Q}$-algebra homomorphism 
${\cal H} \rightarrow {\rm End}_{\mathbb Q}(V)$ sending 
$\widetilde{T}(m_{1}, m_{2})$ to $T(m_{1}, m_{2})$. 
\end{it}
\vspace{2ex}

{\it Proof.} 
The degree of $\widetilde{T}(m_{1}, m_{2})$ in ${\cal H}$ is the number of 
subgroups $G \subset {\cal E}_{\zeta} \otimes \overline{F}_{\zeta}$ 
of type $(m_{1}, m_{2})$, 
and hence by Proposition 3.4, 
the map sending $\widetilde{T}(m_{1}, m_{2})$ to $T(m_{1}, m_{2})$ 
is compatible with the multiplications on ${\cal H}$ 
and ${\rm End}_{\mathbb Q}(V)$. 
\ $\square$ 
\vspace{2ex}

We call the images of this homomorphism 
${\cal H} \rightarrow {\rm End}_{\mathbb Q}(V)$ 
{\it Hecke operators} in our theory. 
\vspace{2ex}

{\sc Theorem 4.2.} 
\begin{it} 
Any Hecke operator is commutative with 
the Galois action on $V \otimes {\mathbb Q}_{l}$, 
and preserves the weight and Hodge filtrations of $V \otimes {\mathbb C}$. 
\end{it}
\vspace{2ex}

{\it Proof.} 
The former assertion follows from Theorem 3.5 (2), 
and we will show the latter one. 
Recall that $T(m_{1}, m_{2})$ on 
$A_{N} \left( \overline{M}_{\zeta}({\mathbb C}); c \right)$ 
is the sum of $\tau_{\zeta, G}$ 
for $G \cong {\mathbb Z}/m_{1} \oplus {\mathbb Z}/m_{2}$. 
Take $z \in \overline{X}_{\zeta, G}({\mathbb Q})$ such that 
$\phi_{\zeta, G}(z) = c \in C_{\zeta}$. 
Then $\phi_{\zeta, G}, \psi_{\zeta, G}$ are holomorphic maps, 
$A_{N} \left( \phi_{\zeta, G} \right), A_{N} \left( \psi_{\zeta, G} \right)$ 
preserve the weight and Hodge filtrations. 
By Proposition 3.1, 
$A_{N} \left( \phi_{\zeta, G} \right)$ is surjective, 
and hence 
$$
A_{N} \left( \phi_{\zeta, G} \right) : 
A_{N} \left( \overline{X}_{\zeta, G}({\mathbb C}); z \right) \rightarrow 
A_{N} \left( \overline{M}_{\zeta}({\mathbb C}); \phi_{\zeta, G}(z) \right)
$$
gives a surjection between each pair of the Hodge components. 
Since 
$$
\tau_{\zeta, G} \circ A_{N} \left( \phi_{\zeta, G} \right) = 
A_{N} \left( \psi_{\zeta, G} \right), 
$$
$\tau_{\zeta, G}$ preserves the Hodge decompositions, 
and hence their sum $T(m_{1}, m_{2})$ has the same property. \ $\square$ 
\vspace{2ex}

The following corollary gives a solution to a problem posed by Manin [M2, 3.3] 
for cusp forms of weight $2$. 
\vspace{2ex}

{\sc Corollary 4.3.} 
\begin{it}
Let $B_{N}$ be as in 2.2 which is the dual space of $V \otimes {\mathbb C}$ 
with action of ${\cal H}$. 
Then the subspace of $B_{N}$ spanned by iterated Shimura integrals 
$$
\int \omega_{1} \cdots \omega_{N} \ 
\left( \mbox{$\omega_{i} :$ holomorphic $1$-forms on 
$\overline{M}_{\zeta}({\mathbb C})$} \right) 
$$
is stable under the action of ${\cal H}$. 
\end{it}
\vspace{2ex}

{\it Proof.} 
Iterated integrals of holomorphic $1$-forms are homotopy functional, 
and by Theorem 4.2, 
the action of ${\cal H}$ preserves the Hodge filtration of $B_{N}$, 
especially its subspace $F^{N}(B_{N})$. \ $\square$ 
\vspace{2ex}

4.2. \ 
Since $V$ is finite dimensional over ${\mathbb Q}$, 
the eigenvalues of each Hecke operator are in $\overline{\mathbb Q}$. 
For each representation 
$\varepsilon : {\cal H} \rightarrow \overline{\mathbb Q}$, 
i.e., ${\mathbb Q}$-algebra homomorphism, 
we define the {\it Hecke component} for $\varepsilon$ as the subspace of 
$V_{\overline{\mathbb Q}} = V \otimes \overline{\mathbb Q}$ given by 
$$
V_{\overline{\mathbb Q}} (\varepsilon) = \bigcap_{h \in {\cal H}} 
\left( \bigcup_{m \in {\mathbb N}} 
{\rm Ker} \left( (h - \varepsilon(h))^{m} \right) \right). 
$$
Then by the commutativity of ${\cal H}$, 
$V_{\overline{\mathbb Q}}$ is decomposed to the direct sum 
of the Hecke components: 
$$
V_{\overline{\mathbb Q}} = \bigoplus_{\varepsilon} 
V_{\overline{\mathbb Q}} (\varepsilon), 
$$
and hence by Theorems 4.2 and 3.5 (1), we have: 
\vspace{2ex}

{\sc Theorem 4.4.} 
\begin{it} 
\begin{itemize}

\item[{\rm (1)}] 
$V_{\overline{\mathbb Q}} (\varepsilon) \otimes \overline{\mathbb Q}_{l}$ 
has Galois action by $\sigma \otimes {\rm id}_{\overline{\mathbb Q}_{l}}$ 
$\left( \sigma \in {\rm Gal} \left( \overline{\mathbb Q}/{\mathbb Q} 
\right) \right)$ on 
$\left( V \otimes {\mathbb Q}_{l} \right) \otimes \overline{\mathbb Q}_{l}$, 
and $V_{\overline{\mathbb Q}} (\varepsilon) \otimes {\mathbb C}$ 
has mixed Hodge structure which are compatible with 
the motivic structure on $V$. 

\item[{\rm (2)}] 
For a prime $l \neq p$, 
the congruence relation 
$$
T(p) = F_{p} + p \left( F_{p}^{-1} \circ T(p,p) \right) 
$$
holds on 
$V_{\overline{\mathbb Q}}(\varepsilon) \otimes \overline{\mathbb Q}_{l}$. 

\end{itemize}
\end{it}

4.3. \ We study the mixedness of $V$ 
in the category of motives over $\overline{\mathbb Q}$, 
where morphisms are considered as 
$\overline{\mathbb Q}$-linear homomorphisms 
compatible with Galois action and with weight and Hodge filtrations. 
For each $\zeta \in P_{n}$ and $c \in C_{\zeta}$, 
let $I(\zeta, c)$ be the augmentation ideal of 
${\mathbb Q} \left[ \pi_{1} \left( \overline{M}_{\zeta}({\mathbb C}); 
c \right) \right]$, 
and put 
$$
I^{m} = \bigoplus_{\zeta, c} I(\zeta, c)^{m} 
$$
which make a decreasing sequence of ${\mathbb Q}$-subspaces of $V$. 
By construction, the action of ${\cal H}$ preserves this filtration. 
\vspace{2ex}

{\sc Theorem 4.5.} 
\begin{it}
For all integers $1 < l < m \leq N$, 
the exact sequence 
$$
0 \rightarrow  \left( I^{l} / I^{m} \right) \otimes \overline{\mathbb Q} 
\rightarrow \left( I / I^{m} \right) \otimes \overline{\mathbb Q} 
\rightarrow \left( I / I^{l} \right) \otimes \overline{\mathbb Q} 
\rightarrow 0 
$$
splits as motives over $\overline{\mathbb Q}$. 
\end{it}
\vspace{2ex}

{\it Proof.} 
Let $p$ be a prime not dividing $n$. 
Then by the famous result of Weil, 
$$
I(\zeta, c) / I(\zeta, c)^{2} \cong 
H_{1} \left( \overline{M}_{\zeta}({\mathbb C}); {\mathbb Q} \right) 
$$
has pure weight $-1$, 
i.e., after scalar-extended to ${\mathbb Q}_{l}$, 
a $p^{d}$th power Frobenius homomorphism 
has eigenvalues with absolute value $p^{d/2}$. 
Hence by the surjectivity of the natural ${\mathbb Q}$-linear homomorphism 
$$
\left( I(\zeta, c) / I(\zeta, c)^{2} \right)^{\otimes m} \rightarrow 
I(\zeta, c)^{m} / I(\zeta, c)^{m+1}, 
$$
each $I^{m} / I^{m+1}$ has pure weight $-m$. 
As seen above, 
$I \otimes \overline{\mathbb Q}$ is the direct sum of 
$$
I_{\overline{\mathbb Q}}(\varepsilon) = V_{\overline{\mathbb Q}}(\varepsilon) 
\cap \left( I \otimes \overline{\mathbb Q} \right) 
$$
as a motive over $\overline{\mathbb Q}$. 
Therefore, to show the assertion, 
it is enough to prove that 
$I_{\overline{\mathbb Q}}(\varepsilon) \cap 
\left( I^{m} \otimes \overline{\mathbb Q} \right)$ 
becomes either $\{ 0 \}$ or $I_{\overline{\mathbb Q}}(\varepsilon)$ 
for any $\varepsilon$ and $m \geq 2$ 
since in this case, 
$\left( I^{l}/I^{m} \right) \otimes \overline{\mathbb Q}$ 
is isomorphic to the direct sum of $I_{\overline{\mathbb Q}}(\varepsilon)$ 
which are contained in $I^{l} \otimes \overline{\mathbb Q}$ 
and not contained in $I^{m} \otimes \overline{\mathbb Q}$. 
Assume on the contrary that 
$$
\{ 0 \} \subsetneq 
W_{1} = I_{\overline{\mathbb Q}}(\varepsilon) \cap 
(I^{m} \otimes \overline{\mathbb Q}) \subsetneq 
W_{2} = I_{\overline{\mathbb Q}}(\varepsilon) 
$$
for some $\varepsilon$ and $m \geq 2$. 
Let $e_{1}$ and $e_{2}$ be the eigenvalues of the Frobenius $F_{p}$ 
on $W_{1} \otimes \overline{\mathbb Q}_{l}$ and on 
$\left( W_{2} / W_{1} \right) \otimes \overline{\mathbb Q}_{l}$ respectively. 
Then by Theorem 4.4 (2), 
$e_{1}$ and $e_{2}$ are the roots of 
$$
x^{2} - \varepsilon(T(p)) x + p \cdot \varepsilon(T(p, p)) = 0
$$ 
such that $p^{1/2} \leq | e_{2} | < p^{m/2} \leq | e_{1} |$. 
Therefore, $| p \cdot \varepsilon(T(p, p)) | = | e_{1} e_{2} | > p$ 
which contradicts with that $T(p, p)$ is an automorphism of $V$ 
with finite order. 
This completes the proof. \ $\square$ 
\vspace{2ex} 

{\sc Theorem 4.6.} 
\begin{it}
Fix $\zeta \in P_{n}$ and $c \in C_{\zeta}$. 
Then for any 
$\gamma \in \pi_{1} \left( \overline{M}_{\zeta}({\mathbb C}); c \right)$ 
and all holomorphic $1$-forms $\omega_{1},...,\omega_{m}$ on 
$\overline{M}_{\zeta}({\mathbb C})$, 
there exist a positive integer $l$, $d_{j} \in \overline{\mathbb Q}$ and 
$\delta_{ij} \in H_{1} \left( \overline{M}_{\zeta}({\mathbb C}); 
{\mathbb Z} \right)$ 
$(1 \leq i \leq m, 1 \leq j \leq l)$ such that 
$$
\int_{\gamma} \omega_{1} \cdots \omega_{m} = 
\sum_{j=1}^{l} d_{j} 
\left( \prod_{i=1}^{m} \int_{\delta_{ij}} \omega_{i} \right). 
$$
\end{it}
\vspace{-1ex} 

{\it Proof.} 
We use results in [H2, \S 6] with some extension. 
Take $\left\{ g(\zeta, c)_{i} \right\}_{i} \subset I(\zeta, c)$ 
giving a basis of $I/I^{m}$ which also provides a section of 
the projection $I/I^{m+1} \rightarrow I/I^{m}$. 
For each element $w_{1} \otimes \cdots \otimes w_{m}$ 
$\left( \mbox{$w_{i}$ : closed $1$-forms on 
$\overline{M}_{\zeta}({\mathbb C})$} \right)$ 
of a basis of 
$$
{\rm Hom} \left( I(\zeta, c)^{m} / I(\zeta, c)^{m+1}, 
{\mathbb Q} \right) \otimes {\mathbb C} \hookrightarrow 
H_{1} \left( \overline{M}_{\zeta}({\mathbb C}); 
{\mathbb C} \right)^{\otimes m} 
$$ 
with Hodge filtration, 
take a set $\left\{ u_{j1},...,u_{jn_{j}} \right\}_{j}$ $(n_{j} < m)$ 
of $1$-forms on $\overline{M}_{\zeta}({\mathbb C})$ such that  
$$
\int w_{1} \cdots w_{m} + \sum_{j} \int u_{j1} \cdots u_{jn_{j}} 
$$
is a homotopy functional iterated integral. 
Then the bilinear form 
\begin{eqnarray*}
\lefteqn{\langle \gamma - 1, w_{1} \otimes \cdots \otimes w_{m} \rangle} 
\\
& = &
\left\{ \begin{array}{l} 
{\displaystyle \int w_{1} \cdots w_{m} + 
\sum_{j} \int u_{j1} \cdots u_{jn_{j}}} 
\ 
\left( \gamma \in \pi_{1} \left( \overline{M}_{\zeta}({\mathbb C}); 
c \right) \right), 
\\
0 
\ 
\left( \gamma \in \pi_{1} \left( \overline{M}_{\zeta'}({\mathbb C}); 
c' \right) \ \mbox{for} \ (\zeta', c') \neq (\zeta, c) 
\right)
\end{array} \right. 
\end{eqnarray*}
gives a retraction of the inclusion 
$\left( I^{m}/I^{m+1} \right) \otimes {\mathbb C} \rightarrow  
\left( I/I^{m+1} \right) \otimes {\mathbb C}$ 
which preserves their Hodge filtrations. 
Therefore, 
the extension data of  
$$
0 \rightarrow I^{m}/I^{m+1} \rightarrow I/I^{m+1} \rightarrow I/I^{m} 
\rightarrow 0
$$
as their mixed Hodge structures over $\overline{\mathbb Q}$ 
is given as an element of 
$$
\frac{{\rm Hom} \left( I/I^{m}, I^{m}/I^{m+1} \right) \otimes {\mathbb C}} 
{F^{0} \left( {\rm Hom} \left( I/I^{m}, 
I^{m}/I^{m+1} \right) \otimes {\mathbb C} \right) 
+ 
{\rm Hom} \left( I/I^{m}, I^{m}/I^{m+1} \right) \otimes \overline{\mathbb Q}} 
$$ 
which sends $g(\zeta, c)_{i}$ to 
$$
w_{1} \otimes \cdots \otimes w_{m} \mapsto 
\langle g(\zeta, \alpha)_{i}, w_{1} \otimes \cdots \otimes w_{m} \rangle. 
$$
Since $\omega_{1},...,\omega_{m}$ are holomorphic, 
${\displaystyle \int \omega_{1} \cdots \omega_{m}}$ is homotopy functional, 
and any element of 
$F^{0} \left( {\rm Hom} \left( I/I^{m}, 
I^{m}/I^{m+1} \right) \otimes {\mathbb C} \right)$ 
sends each element of $I/I^{m}$ to 
$$
F^{m} \left( {\rm Hom} \left( I^{m}/I^{m+1}, {\mathbb Q} \right) \otimes 
{\mathbb C} \right) \ni \omega_{1} \otimes \cdots \otimes \omega_{m} 
\mapsto 0 
$$
because 
$F^{m} \left( {\rm Hom} \left( I/I^{m}, {\mathbb Q} \right) \otimes 
{\mathbb C} \right) = \{ 0 \}$. 
Therefore, by Theorem 4.5, 
$$
\omega_{1} \otimes \cdots \otimes \omega_{m} \mapsto 
\langle g(\zeta, c)_{i}, \omega_{1} \otimes \cdots \otimes \omega_{m} 
\rangle
$$
belongs to $\left( I^{m}/I^{m+1} \right) \otimes \overline{\mathbb Q}$, 
and hence becomes a $\overline{\mathbb Q}$-linear sum of 
${\displaystyle \prod_{i=1}^{m} \int_{\delta_{i}} \omega_{i}}$ 
for some 
$\delta_{i} \in H_{1} \left( \overline{M}_{\zeta}({\mathbb C}); 
{\mathbb Z} \right)$. \ $\square$ 
\vspace{2ex} 

Theorem 4.6 and calculation of iterated Shimura integrals in [M1] 
imply the following: 
\vspace{2ex} 

{\sc Corollary 4.7.} 
\begin{it}
For cusp forms 
$$
\omega_{i} = \sum_{l=1}^{\infty} c_{i}(l) e^{2 \pi \sqrt{-1} l \tau/n} \ 
(1 \leq i \leq m) 
$$
of weight $2$ and level $n$, 
and for $a \in \Gamma(n) \left( \sqrt{-1} \cdot \infty \right)$, 
the multiple $L$-value 
$$
\sum_{0 < l_{1} < \cdots < l_{m}} 
\frac{c_{1}(l_{m}-l_{m-1}) \cdots c_{m}(l_{1})}{l_{m} \cdots l_{1}} 
e^{2 \pi \sqrt{-1} l_{m} a / n} 
$$
becomes a $\overline{\mathbb Q}$-linear sum of the products 
$$
\prod_{i=1}^{m} \left( \sum_{l_{i} = 1}^{\infty} 
\frac{c_{i}(l_{i})}{l_{i}} e^{2 \pi \sqrt{-1} l_{i} a_{i} / n} \right) 
$$ 
of $L$-values for some 
$a_{i} \in \Gamma(n) \left( \sqrt{-1} \cdot \infty \right)$. 
\end{it}
\vspace{2ex} 

{\it Remark.} 
As is seen in the above consideration, 
``$\overline{\mathbb Q}$'' in Theorems 4.5, 4.6 and Corollary 4.7 
can be replaced with a finite extension of ${\mathbb Q}$ over which 
the decomposition $I = \bigoplus_{\varepsilon} I(\varepsilon)$ is given. 
\vspace{2ex}

\begin{center}
{\bf 5. Algebraicity of periods of divisors}
\end{center}

5.1. \ 
Let $\overline{M}_{\zeta}$ be the above compact modular curve of 
level $n \geq 3$ and genus $g > 0$ which is defined over ${\mathbb Q}(\zeta)$, 
where $\zeta$ denotes a primitive $n$th root of $1$. 
Denote by $J$ the jacobian variety of $\overline{M}_{\zeta}$ on which $-1$ 
acts as $x \mapsto -x$, 
and put $X^{-} = (-1)^{*}(X)$ for cycles $X$ on $J$. 
Fix a cusp $c$ on $\overline{M}_{\zeta}$, 
and let $\iota$ be the embedding $\overline{M}_{\zeta} \hookrightarrow J$ 
sending $x$ to $[x - c]$. 
Then the aim of this section is to show the following: 
\vspace{2ex}

{\sc Theorem 5.1.} 
\begin{it} 
For a divisor ${\cal D}$ on $J$, 
let $D$ be a divisor on $\overline{M}_{\zeta}$ such that 
${\cal O}_{\overline{M}_{\zeta}}(D) \cong \iota^{*} 
\left( {\cal O}_{J}({\cal D}) \otimes {\cal O}_{J}({\cal D}^{-}) \right)$, 
and put $D_{0} = D - \deg(D) \cdot c$. 
Then for any holomorphic $1$-form $\omega$ on $\overline{M}_{\zeta}$, 
there exist algebraic numbers $a_{i}$ such that 
$$
\int_{D_{0}} \omega = \sum_{i} a_{i} \cdot \int_{\gamma_{i}} \omega, 
$$
where $\gamma_{i}$ are basis of 
$H_{1} \left( \overline{M}_{\zeta}({\mathbb C}); {\mathbb Q} \right)$. 
\end{it} 
\vspace{2ex}

This theorem is obtained by combining a result of Pulte [P] 
on Ceresa cycles [C] with results in Section 4. 
Then the above $a_{i}$ are computed by describing Hecke operators 
on these motives, and by Abel's theorem, 
if one will know that there is a $a_{i} \not\in {\mathbb Q}$, 
then $D_{0}$ has infinite order in $J$. 
We hope that the computation of $a_{i}$ is applied to study 
the nontriviality of the modular Ceresa cycle 
$\overline{M}_{\zeta} - \left( \overline{M}_{\zeta} \right)^{-}$ 
(modulo algebraic equivalence) and of $D_{0}$ (modulo torsion). 
\vspace{2ex}

5.2. \ 
We give a proof of Theorem 5.1. 
Let 
$I(\zeta, c) \subset {\mathbb Q} 
\left[ \pi_{1} \left( \overline{M}_{\zeta}({\mathbb C}); c \right) \right]$ 
be as above. 
Then 
$V_{i} = H_{i} \left( \overline{M}_{\zeta}({\mathbb C}); 
{\mathbb Q} \right)$ 
satisfies that 
$$
I(\zeta, c) / I(\zeta, c)^{2} \cong V_{1} \cong V_{1}^{\vee}(1), \ 
I(\zeta, c)^{2} / I(\zeta, c)^{3} \cong 
\left. \left( \wedge^{2} V_{1} \right) \right/ V_{2}, 
$$
where $V_{1}^{\vee}(1)$ denotes the Tate twist of the dual space of $V_{1}$, 
and $V_{2} \rightarrow \wedge^{2} V_{1}$ is the dual map of the cup product. 
Since the modular Ceresa cycle 
$\overline{M}_{\zeta} - \left( \overline{M}_{\zeta} \right)^{-}$ 
is homologically trivial, 
it gives rise to an element ${\cal E}$ of 
$$
{\rm Ext} \left( {\mathbb Q}, 
H^{2g-3} \left( J({\mathbb C}), {\mathbb Q} \right)(g-1) \right) 
\cong 
{\rm Ext} \left( {\mathbb Q}, \wedge^{3} V_{1}(-1) \right) 
$$
in the abelian category of mixed ${\mathbb Q}$-Hodge structures. 
Furthermore, by a result of Pulte [P] (see also [H2, 3]), 
${\cal E}$ becomes the twice of the elements corresponding to 
the natural exact sequence 
$$
0 \rightarrow I(\zeta, c)^{2} / I(\zeta, c)^{3} \rightarrow 
I(\zeta, c) / I(\zeta, c)^{3} \rightarrow I(\zeta, c) / I(\zeta, c)^{2} 
\rightarrow 0
$$
under an inclusion 
$$
\wedge^{3} V_{1} \rightarrow 
V_{1} \otimes 
\left( \left. \left( \wedge^{2} V_{1} \right) \right/ V_{2} \right) 
\cong 
\left( I(\zeta, c) / I(\zeta, c)^{2} \right)^{\vee}(1) \otimes 
\left( I(\zeta, c)^{2} / I(\zeta, c)^{3} \right) 
$$
sending $x \wedge y \wedge z$ to 
$x \otimes (y \wedge z) + y \otimes (z \wedge x) + z \otimes (x \wedge y)$. 
We consider the category of mixed $\overline{\mathbb Q}$-Hodge structures 
as that of finite dimensional $\overline{\mathbb Q}$-vector space with 
weight filtration $W_{*}$ and Hodge filtration $F^{*}$ over ${\mathbb C}$. 
Note that in this category, a morphism $f : V_{1} \rightarrow V_{2}$ 
is required to be a $\overline{\mathbb Q}$-linear map strongly compatible with 
weight and Hodge filtrations, i.e., 
$$
f \left( W_{n}(V_{1}) \right) = f(V_{1}) \cap W_{n}(V_{2}), \ 
f \left( F^{p}(V_{1} \otimes {\mathbb C}) \right) = 
f(V_{1} \otimes {\mathbb C}) \cap F^{p}(V_{2} \otimes {\mathbb C}). 
$$
Then mixed $\overline{\mathbb Q}$-Hodge structures make an abelian category. 
By the splitting constructed in Theorem 4.5, 
the natural injection 
$I(\zeta, c) / I(\zeta, c)^{2} \rightarrow I / I^{2}$ and 
the natural surjection 
$I / I^{3} \rightarrow I(\zeta, c) / I(\zeta, c)^{3}$ 
give rise to a splitting of the above exact sequence as motives 
over $\overline{\mathbb Q}$. 
Therefore, ${\cal E}$ is trivial in the abelian category of 
mixed $\overline{\mathbb Q}$-Hodge structures 
since $\wedge^{3} V_{1}$ becomes a direct summand of 
$\left( I(\zeta, c) / I(\zeta, c)^{2} \right)^{\vee}(1) \otimes 
\left( I(\zeta, c)^{2} / I(\zeta, c)^{3} \right)$. 
Under the natural map 
$$
{\rm Ext} \left( {\mathbb Q}, \wedge^{3} V_{1}(-1) \right) \times 
H^{2}(J, {\mathbb Q})(1) \rightarrow 
{\rm Ext} \left( {\mathbb Q}, V_{1} \right), 
$$
the product $E$ of ${\cal E}$ and $[{\cal D}]$ corresponds to the divisor 
$$
\left( \overline{M}_{\zeta} - \left( \overline{M}_{\zeta} \right)^{-} \right) 
\cdot \left( {\cal D} + {\cal D}^{-} \right) = 
2 \overline{M}_{\zeta} \cdot \left( {\cal D} + {\cal D}^{-} \right) = 2 D.
$$
Therefore, $E$ is trivial in the abelian category of 
mixed $\overline{\mathbb Q}$-Hodge structures. 
In the expression 
$$
{\rm Ext} \left( {\mathbb Q}, V_{1} \right) \cong 
\frac{V_{1} \otimes {\mathbb C}}
{F^{0} \left( V_{1} \otimes {\mathbb C} \right) + V_{1}}, 
$$
$E$ corresponds to the integration on $D_{0}$, 
and hence the triviality of $E$ over $\overline{\mathbb Q}$ implies 
the assertion 
since any element of $F^{0} \left( V_{1} \otimes {\mathbb C} \right)$ 
maps $\omega$ to $0$ if $\omega$ is holomorphic. 
\ $\square$ 
\vspace{2ex}

\begin{center} 
{\sc Acknowledgments}
\end{center}

The author would like to deeply thank Akio Tamagawa for kindly pointing out 
mistakes on the construction of Hecke action in the previous version. 
\vspace{4ex}

{\sc Department of Mathematics,
Graduate School of Science and Engineering,
Saga University, Saga 840-8502, Japan} 

{\it E-mail address}: ichikawa@ms.saga-u.ac.jp 
\vspace{2ex}

\begin{center} 
{\sc References}
\end{center}

\begin{itemize} 

\item[{[C]}] 
{\sc G. Ceresa,} 
$C$ is not algebraically equivalent to $C^{-}$ in its Jacobian, 
{\it Ann. of Math.} {\bf 117} (1983), 285--291. 

\item[{[Ch1]}] 
{\sc K. T. Chen,} 
Iterated path integrals, 
{\it Bull. Amer. Math. Soc.} {\bf 83} (1977), 831--879. 

\item[{[Ch2]}] 
{\sc K. T. Chen,} 
Extension of $C^{\infty}$ functions algebra by integrals 
and Malcev completion of $\pi_{1}$, 
{\it Adv. in Math.} {\bf 23} (1977), 181--210. 

\item[{[D1]}] 
{\sc P. Deligne,} 
Formes modulaires et repr\'{e}sentations $l$-adiques, 
{\it Expos\'{e} 355, S\'{e}minaire N. Bourbaki, 1968/69, 
Lecture Notes in Math.} 
{\bf 179}, Springer-Verlag, 1969, pp. 139--172. 

\item[{[D2]}] 
{\sc P. Deligne,} 
Le groupe fondamental de la droite projective moins trois points, 
in: Y. Ihara, K. Ribet and J. P. Serre, (eds.), 
{\it Galois groups over ${\mathbb Q},$ Publ. MSRI} {\bf 16} 1989, pp. 79--298. 

\item[{[DR]}] 
{\sc P. Deligne} and {\sc M. Rapoport,}
Les sch\'{e}mas de modules de courbes elliptiques, 
{\it Modular Functions of One Variable II, 
Lecture Notes in Math.} 
{\bf 349}, Springer-Verlag, 1973, pp. 143--316. 

\item[{[G]}] 
{\sc A. Grothendieck} et al., 
{\it S\'{e}minaire de G\'{e}om\'{e}trie Alg\'{e}brique, SGA 1. 
Rev\^{e}tements \'{e}tale et Groupe Fondemental, 
Lecture Notes in Math.} {\bf 224}, Springer-Verlag, 1971. 

\item[{[H1]}] 
{\sc R. Hain,} 
The de Rham homotopy of complex algebraic varieties I, II, 
{\it K-theory} {\bf 1} (1987), 271--324, 481--497. 

\item[{[H2]}] 
{\sc R. Hain,} 
The geometry of the mixed Hodge structure on the fundamental group, 
{\it Proc. Symp. Pure Math.} {\bf 6-2} (1987), 247--282. 

\item[{[H3]}] 
{\sc R. Hain,} 
Completions of mapping class groups and the cycle $C - C^{-}$, 
{\it Contemp. Math.} {\bf 150} (1993), 75--105. 

\item[{[I]}] 
{\sc J. Igusa,} 
Kroneckerian model of fields of elliptic modular functions, 
{\it Amer. J. Math.} {\bf 81} (1959), 561--577. 

\item[{[KM]}] 
{\sc N. M. Katz} and {\sc B. Mazur,} 
{\it Arithmetic moduli of elliptic curves, Ann. Math. Stud.} 
{\bf 108} Princeton University Press, 1985. 

\item[{[M1]}] 
{\sc Y. I. Manin,} 
Iterated integrals of modular forms and noncommutative modular symbols, 
in: V. Ginzburg, (ed.), 
{\it Algebraic Geometry and Number Theory, Progr. Math.} 
{\bf 253}, Birkh\"{a}user, 2006, pp. 565--597. 

\item[{[M2]}] 
{\sc Y. I. Manin,} 
Iterated Shimura integrals, 
arXiv:math/0507438v1. 

\item[{[Mo]}] 
{\sc J. Morgan,} 
The algebraic topology on smooth algebraic varieties, 
{\it Publ. IHES} {\bf 48} (1978), 137--204. 

\item[{[P]}] 
{\sc M. Pulte,} 
The fundamental group of a Riemann surface: 
mixed Hodge structures and algebraic cycles, 
{\it Duke Math. J.} {\bf 57} (1988), 721--760. 

\item[{[S]}] 
{\sc A. Shiho,} 
Crystalline fundamental groups. I. 
Isocrystals on log crystalline site and log convergent site, 
{\it J. Math. Soc. Univ. Tokyo} {\bf 7} (2000), 509--656. 

\item[{[Sm]}] 
{\sc G. Shimura,} 
{\it Introduction to the arithmetic theory of automorphic functions,} 
Iwanami Shoten Publishers and Princeton University Press, 1971. 

\end{itemize}

\end{document}